\documentclass[a4paper,10pt,
colorlinks,urlcolor=black,linkcolor=black,citecolor=black,
]{article}
\usepackage{amsmath,amssymb}
\usepackage[dvips]{color}
\usepackage{cite}
\usepackage{hangcaption}
\usepackage{amsmath,amssymb}
\usepackage[dvips]{color}
\usepackage{cite}
\usepackage{eepic}
\usepackage{epic}
\usepackage{hangcaption}
\usepackage{amsmath}
\usepackage{amssymb}
\usepackage{amscd}
\usepackage{amsthm}
\usepackage{multicol}
\usepackage[dvipdfmx,%
bookmarks=true,%
bookmarksnumbered=true,%
setpagesize=false,%
pdfkeywords={TeX; dvipdfmx; hyperref; color;}]{hyperref}
\def\cases{\left\{\begin{array}{ll}}
\def\endcases{\end{array}\right.}

\def\bigtimes{\mathop{\mbox{\Large $\times$}}}
\begin{document}
\setcounter{page}{1}
\vskip1.5cm
\begin{center}
{\Large \bf
ANOVA (analysis of variance) 
in the quantum linguistic formulation of statistics
}
\vskip0.5cm
{\rm
\large
Shiro Ishikawa
}
\\
\vskip0.2cm
\rm
\it
Department of Mathematics, Faculty of Science and Technology,
Keio University,
\\ 
3-14-1, Hiyoshi, Kouhoku-ku Yokohama, Japan.
E-mail:
ishikawa@math.keio.ac.jp
\end{center}
\par
\rm
\vskip0.3cm
\par
\noindent
{\bf Abstract}
\normalsize
\vskip0.3cm
\par
\noindent
Recently, we proposed quantum language (or, measurement theory), which is characterized as the linguistic turn of the Copenhagen interpretation of quantum mechanics. We believe that this language has a great powet of
description,
and therefore, even statistics can be described by quantum language.
In this paper, we show that
ANOVA ( analysis of variance (one-way and two-way)) can be formulated in quantum language.
Since quantum language is suited for theoretical arguments, we believe that our results are visible and understandable.
For example,
we can answer the question
"What kind of role does Kolmogorov's probability theory play in ANOVA?"
That is,  the readers find that Kolmogorov's probability theory is merely used in order to calculate multi-dimenstional Gauss integrals, and thus, they can avoid to confuse the relation between Kolmogorov's probability theory and statistics.

\par
\noindent
(Key words: 
Quantum language, Statistical hypothesis testing, ANOVA, F-distribution,
Student's t-distribution, Chi-squared distribution, 
)

\vskip1.0cm

\par

\def\Cal{\cal}
\def\bigstimes{\text{\large $\: \boxtimes \,$}}

\par
\noindent

\vskip0.2cm
\par
\noindent
\par
\noindent
\section{
Introduction
}
%

\rm
\par
\par
\noindent

\par
\noindent
\subsection{
Quantum language
(Axioms
and
Interpretation)
}
As mentioned in the above abstract, our purpose is to understand
ANOVA
( analysis of variance )
in terms of quantum language, which is proposed in \cite{Ishi2}-\cite{Ishi9}.
\par
According to ref.\cite{Ishi9},
we shall mention the overview of quantum language
(or, measurement theory, in short, MT).
\par
\par
\rm
Quantum language is characterized as the linguistic turn of the Copenhagen interpretation of quantum mechanics({\it cf.} refs. \cite{Ishi5},
{{{}}}{\cite{Neum}}).
Quantum language (or, measurement theory ) has two simple rules
(i.e. Axiom 1(concerning measurement) and Axiom 2(concerning causal relation))
and the linguistic interpretation (= how to use the Axioms 1 and 2). 
That is,
\begin{align}
\underset{\mbox{(=MT(measurement theory))}}{\fbox{Quantum language}}
=
\underset{\mbox{(measurement)}}{\fbox{Axiom 1}}
+
\underset{\mbox{(causality)}}{\fbox{Axiom 2}}
+
\underset{\mbox{(how to use Axioms)}}{\fbox{linguistic interpretation}}
\label{eq1}
\end{align}
({\it cf.} refs.
{{{}}}{\cite{Ishi2}-\cite{Ishi9}}).
\par
This theory is formulated in a certain $C^*$-algebra ${\cal A}$({\it cf.} ref.
{{{}}}{\cite{Saka}}), and is classified as follows:
\begin{itemize}
\item[(A)]
$
\quad
\underset{\text{\scriptsize }}{\text{MT}}
$
$\left\{\begin{array}{ll}
\text{quantum MT$\quad$(when ${\cal A}$ is non-commutative)}
\\
\\
\text{classical MT
$\quad$
(when ${\cal A}$ is commutative, i.e., ${\cal A}=C_0(\Omega)$)}
\end{array}\right.
$
\end{itemize}
where $C_0(\Omega)$
is
the $C^*$-algebra composed of all continuous 
complex-valued functions vanishing at infinity
on a locally compact Hausdorff space $\Omega$.

Since our concern in this paper is concentrated to 
the usual statistical hypothesis test
methods in statistics,
we devote ourselves to the commutative $C^*$-algebra $C_0(\Omega)$,
which is quite elementary.
Therefore, we believe that all statisticians
can understand our assertion
(i.e.,
a new viewpoint of the confidence interval
methods
).

Let $\Omega$ is a locally compact Hausdorff space, which is also called
a state space. And thus, an element $\omega (\in \Omega )$ is said to be a state.
Let $C(\Omega)$ be the $C^*$-algebra composed of all bounded continuous 
complex-valued functions on a locally compact Hausdorff space $\Omega$.
The norm $\| \cdot \|_{C(\Omega )}$ is usual, i.e.,
$\| f \|_{C(\Omega )} = \sup_{\omega \in \Omega } |f(\omega )|$
$(\forall f \in C(\Omega ))$.
\par
\rm
Motivated by Davies' idea ({\it cf.} ref.
{{{}}}{\cite{Davi}}) in quantum mechanics,
an observable ${\mathsf O}=(X, {\mathcal F}, F)$ in $C_0(\Omega )$
(or, precisely, in $C(\Omega )$) is defined as follows:
\begin{itemize}
\item[(B$_1$)]
$X$ is a topological space. 
${\mathcal F} ( \subseteq 2^X$(i.e., the power set of $X$) is a field,
that is, it satisfies the following conditions (i)--(iii):
(i):
$\emptyset \in {\cal F}$, 
(ii):$\Xi \in {\mathcal F} \Longrightarrow X\setminus \Xi  \in 
{\mathcal F}$,
(iii):
$\Xi_1, \Xi_2,\ldots, \Xi_n \in {\mathcal F} \Longrightarrow \cup_{k=1}^n \Xi_k \in {\mathcal F}$.
\item[(B$_2$)]
The map $F: {\cal F} \to C(\Omega )$ satisfies that 
\begin{align}
0 \le [F(\Xi )](\omega ) \le 1, \quad [F(X )](\omega )=1
\qquad
(\forall \omega \in \Omega )
\nonumber \end{align}
and moreover, 
if
\begin{align}
\Xi_1, \Xi_2,\ldots, \Xi_n, \ldots \in {\mathcal F},
\quad
\Xi_m \cap \Xi_n = \emptyset \quad( m \not= n ),
\quad
\Xi = \cup_{k=1}^\infty \Xi_k \in {\mathcal F},
\nonumber \end{align} 
then, it holds
\begin{align}
[F(\Xi)](\omega) = \lim_{n \to \infty } \sum_{k=1}^n [F(\Xi_k )](\omega )
\quad
(\forall \omega \in \Omega )
\nonumber \end{align}
\end{itemize}
Note that Hopf extension theorem
({\it cf.}
ref.
{{{}}}{\cite{Yosi}})
guarantees that
$(X, {\cal F}, [F(\cdot)](\omega))$
is regarded as the mathematical probability space.
\par

\vskip0.5cm

%
%

\par
Now we shall briefly explain "quantum language (\ref{eq1})" in classical systems as follows:
A measurement of an observable
${\mathsf O}=(X, {\mathcal F}, F)$
for a system with a state $\omega (\in \Omega )$
is denoted by
${\mathsf M}_{C_0(\Omega)} ({\mathsf O}, S_{[\omega]})$.
By the measurement, a measured value $x (\in X)$ is obtained as follows:
\par
\noindent
\bf
Axiom 1
\rm
(Measurement)
\begin{itemize}
\item{}
\sl
The probability that a measured value $x$
$( \in X)$ obtained by the measurement 
${\mathsf{M}}_{{{C_0(\Omega)}}} ({\mathsf{O}}$
${ \equiv} (X, {\cal F}, F),$
{}{$ S_{[\omega_0]})$}
belongs to a set 
$\Xi (\in {\cal F})$ is given by
$
[F(\Xi) ](\omega_0 )
$.
\end{itemize}
\rm
\par
\noindent
\par
\noindent
\bf
Axiom 2
\rm
(Causality)
\begin{itemize}
\item{}
\sl
The causality is represented by a Markov operator
$\Phi_{21} : C_0(\Omega_2 ) \to C_0(\Omega_1 )$.
Particularly, the deterministic causality
is represented by a continuous map
$\pi_{12} : \Omega_1 \to \Omega_2$
\end{itemize}
\par
\noindent
\bf
Interpretation
\rm
(Linguistic interpretation).
Although there are several linguistic rules in quantum language, the following is the most important:
\begin{itemize}
\item{}
\sl
Only one measurement is permitted.
\end{itemize}
\rm
In order to read this paper,
it suffices to understand the above three
( particularly, Axiom 1).
For the further arguments, see refs.
{{{}}}{\cite{Ishi2}-\cite{Ishi9}}.
\vskip0.5cm

\par
\noindent
\bf
Example 1
\rm
[A kind of normal observable].
Let $n$ be a natural number.
Then, we get a kind of normal observable
${\mathsf O}_G^n = ({\mathbb R}^n, {\mathcal B}_{\mathbb R}^n, {{{G}}^n} )$ 
in $C_0({\mathbb R} \times {\mathbb R}_+)$.
That is,
\par
\noindent
\begin{align}
&
[{{{G}}}^n
(\bigtimes_{k=1}^n \Xi_k)]
({}\omega{})
=
\bigtimes_{k=1}^n
[{{{G}}}(\Xi_k)](\omega)
\nonumber
\\
=
&
\frac{1}{({{\sqrt{2 \pi }\sigma{}}})^n}
\underset{{\bigtimes_{k=1}^n \Xi_k }}{\int \cdots \int}
\exp[{}- \frac{\sum_{k=1}^n ({}{}{x_k} - {}{\mu}  {})^2 
}
{2 \sigma^2}    {}] d {}{x_1} d {}{x_2}\cdots dx_n
\label{eq2}
\\
&
\qquad 
({}\forall  \Xi_k \in {\cal B}_{{\mathbb R}}
\mbox{(=Borel field in ${\mathbb R}$)},
({}k=1,2,\ldots, n),
\quad
\forall   {}{\omega}=(\mu, \sigma )    \in \Omega = {\mathbb R}\times {\mathbb R}_+{}).
\nonumber
\end{align}
Fisher's maximum likelihood method
({\it cf}. refs. {\cite{Ishi3}-\cite{Ishi9}}) urges us to
define the maps
$\overline{\mu}: {\mathbb R}^n \to {\mathbb R}$, 
${\overline{\sigma}}: {\mathbb R}^n \to {\mathbb R}$
and
${\overline{SS}}: {\mathbb R}^n \to {\mathbb R}$
such that
\begin{align}
&
\overline{\mu}
(x) =
\overline{\mu}
(x_1,x_2,\ldots , x_n ) =
\frac{x_1 + x_2 + \cdots + x_n}{n}
\quad( \forall x=(x_1,x_2,\ldots , x_n ) \in {\mathbb R}^n )
\label{eq3}
\\
&
{{\overline{\sigma}}}
(x) =
{{\overline{\sigma}}}
(x_1,x_2,\ldots , x_n ) =
\sqrt{
\frac
{\sum_{k=1}^n ( x_k - 
\overline{\mu}
(x))^2}
{n}
}
\quad( \forall x=(x_1,x_2,\ldots , x_n ) \in {\mathbb R}^n)
\label{eq4}
\intertext{and}
&
{{\overline{SS}}}
(x) =
{{\overline{SS}}}
(x_1,x_2,\ldots , x_n ) =
{\sum_{k=1}^n ( x_k - 
\overline{\mu}
(x))^2}
=n (\overline{\sigma}(x))^2
\quad( \forall x=(x_1,x_2,\ldots , x_n ) \in {\mathbb R}^n)
\label{eq5}
\end{align}
Thus, we have the@following two
image observables
$\overline{\mu}({\mathsf O}_G^n) $
$= ({\mathbb R}, {\mathcal B}_{\mathbb R}, {{{G}}^n} \circ \overline{\mu}^{-1} )$
and
${{\overline{SS}}}({\mathsf O}_G^n) $
$= ({\mathbb R}_+, {\mathcal B}_{{\mathbb R}_+}, {{{G}}^n} \circ {{\overline{SS}}}^{-1} )$
in $C_0({\mathbb R} \times {\mathbb R}_+)$,
which are obatained by the formulas
of
Gauss integrals.
\begin{align}
&
[({{{G}}^n} \circ \overline{\mu}^{-1})(\Xi_1)](\omega)
=
\frac{1}{({{\sqrt{2 \pi }\sigma{}}})^n}
\underset{
\{ x \in {\mathbb R}^n \;:\; {\overline{\mu}}(x) \in \Xi_1 \}}
{\int \cdots \int}
\exp[{}- \frac{\sum_{k=1}^n ({}{}{x_k} - {}{\mu}  {})^2 
}
{2 \sigma^2}    {}] d {}{x_1} d {}{x_2}\cdots dx_n
\nonumber
\\
=
&
\frac{\sqrt{n}}{{\sqrt{2 \pi }\sigma{}}}
\int_{{\Xi_1}} \exp[{}- \frac{n({}{}{x} - {}{\mu}  {})^2 }{2 \sigma^2}    {}] d {}{x}
\label{eq6}
\intertext{and}
&
[({{{G}}^n} \circ {{{\overline{\sigma}}}}^{-1})\big((0,\eta]\big)](\omega)
=
\frac{1}{({{\sqrt{2 \pi }\sigma{}}})^n}
\underset{
0< {\overline{\sigma}}(x) \le \eta }
{\int \cdots \int}
\exp[{}- \frac{\sum_{k=1}^n ({}{}{x_k} - {}{\mu}  {})^2 
}
{2 \sigma^2}    {}] d {}{x_1} d {}{x_2}\cdots dx_n
\nonumber
\\
=
&
\int_0^{n \eta^2/\sigma^2} p^{{\chi}^2}_{n-1}({ x} ) {dx}
\label{eq7}
\\
&
\quad
(
{}\forall  {\Xi_1} \in {\cal B}_{{\mathbb R}}
\mbox{(=Borel field in ${\mathbb R}$)},
\;\;
\forall \eta >
0,
\quad
\forall   {}{\omega} =(\mu, \sigma)   \in \Omega \equiv {\mathbb R}{}\times {\mathbb R}_+).
\nonumber
\end{align}
Here, $p^{{\chi}^2}_{n-1}({ x} )$ is the chi-squared distribution with $n-1$ degrees of freedom. That is,
\begin{align}
p^{{\chi}^2}_{n-1}({ x} )
=
\frac{x^{(n-1)/2-1}e^{-x/2}}{2^{(n-1)/2} \Gamma ((n-1)/2)}
\quad ( x > 0)
\label{eq8}
\end{align}
where $\Gamma$ is the gamma function.

\bf
\par
\noindent
Remark 1
\rm
[The formulas of Gauss integrals].
Although the above (\ref{eq6}) and (\ref{eq7}) can be obtained by direct calculations, we consider that
the calculation in the framework of Kolmogorov's probability theory (ref. \cite{Kolm}) is the most elegant.
This kind of problem
(i.e.,
the formulas of Gauss integrals)
will be repeatedly discussed in this paper
({\it cf}. Remark 2 later).

\par
\vskip1.0cm
\par
\rm
\subsection{The reverse relation between confidence interval and statistical hypothesis testing}
\par
\noindent
\par
Let
${\mathsf O} = ({}X, {\cal F} , F{}){}$
be an observable
formulated in a
commutative $C^*$-algebra
${C_0(\Omega)}$.
Let $X$ be  a topological space.
Let $\Theta$ be a locally compact space with the 
semi-distance $d^x_{\Theta}$
$(\forall x \in X)$,
that is,
for each $x\in X$,
the map
$d^x_{\Theta}: \Theta^2 \to [0,\infty)$
satisfies that
(i):$d^x_\Theta (\theta, \theta )=0$,
(ii):$d^x_\Theta (\theta_1, \theta_2 )$
$=d^x_\Theta (\theta_2, \theta_1 )$,
(ii):$d^x_\Theta (\theta_1, \theta_3 )$
$\le d^x_\Theta (\theta_1, \theta_2 )
+
d^x_\Theta (\theta_2, \theta_3 )
$.

Let
$E:X \to \Theta$
and
$\pi: \Omega \to \Theta$ be
continuous maps,
which are respectively called
an estimator and
a quantity.
Let
$\alpha$
be a real number such that
$0 < \alpha \ll 1$,
for example,
$\alpha = 0.05$.
For any state
$ \omega ({}\in \Omega)$,
define
the positive number
$\eta^\alpha_{\omega}$
$({}> 0)$
such that:
\begin{align}
\eta^\alpha_{\omega}
&
=
\inf
\{
\eta > 0:
[F(\{ x \in X \;:\; 
d^x_\Theta ( E(x) , \pi( \omega ) )
\ge \eta
\}
)](\omega )
\le \alpha
\}
\nonumber
\\
\Big(
&=
\inf
\{
\eta > 0:
[F(\{ x \in X \;:\; 
d^x_\Theta ( E(x) , \pi( \omega ) )
< \eta
\}
)](\omega )
\ge 1- \alpha
\}
\Big)
\label{eq9}
\end{align}
Then Axiom 1 says that:
\rm
\begin{enumerate}
\item[(C$_1$)]
\it
\sl
the probability,
that
the measured value $x$
obtained
by the measurement
${\mathsf M}_{C_0(\Omega)} \big({}{\mathsf O}:= ({}X, {\cal F} , F{})  ,$
$ S_{[\omega_0 {}] } \big)$
satisfies the following
condition (\ref{eq10}),
is more than or equal to
$1-\alpha$
({}e.g., $1-\alpha= 0.95${}).
\begin{align}
d^x_\Theta (E(x),  \pi(\omega_0){}) <  {\eta }^\alpha_{\omega_0}  .
\label{eq10}
\end{align}
\end{enumerate}
or equivalently,
\begin{enumerate}
\item[(C$_2$)]
\it
\sl
the probability,
that
the measured value $x$
obtained
by the measurement
${\mathsf M}_{C_0(\Omega)} \big({}{\mathsf O}:= ({}X, {\cal F} , F{})  ,$
$ S_{[\omega_0 {}] } \big)$
satisfies the following
condition (\ref{eq11}),
is less than or equal to
$\alpha$
({}e.g., $\alpha= 0.05${}).
\begin{align}
d^x_\Theta (E(x),  \pi(\omega_0){}) \ge  {\eta }^\alpha_{\omega_0}  .
\label{eq11}
\end{align}
\end{enumerate}


\rm
\par
\noindent
\bf
Theorem 1
\rm
[Confidence interval and statistical hypothesis testing
({\it cf.} ref. \cite{Ishi9})
].
\sl
Let
${\mathsf O} = ({}X, {\cal F} , F{}){}$
be an observable
formulated in a
commutative $C^*$-algebra
${C_0(\Omega)}$.
Let
$E:X \to \Theta$
and
$\pi: \Omega \to \Theta$ be
an estimator and
a quantity
respectively.
Let $\eta_\omega^\alpha$ be as defined in the formula (\ref{eq9}).
\par
\noindent
From the ($C_1$), we assert "the confidence interval method" as follows:
\begin{enumerate}
\rm
\item[(D$_1$)]
[The confidence interval method].
\sl
\it
\sl
For any $x \in X$, define
\begin{align}
D_{x}^{1- \alpha}
=
\{
\pi(\omega)
(\in
\Theta)
:
d^x_\Theta ({}E(x),
\pi(\omega )
)
<
\eta^{1- \alpha}_{\omega }
\}.
\label{eq12} 
\end{align}
which is called
\it
the $({}1- \alpha{})$-confidence interval. 
\sl
Let
$x (\in X)$
be
a measured value $x$
obtained
by the measurement
${\mathsf M}_{C_0(\Omega)} \big({}{\mathsf O}:= ({}X, {\cal F} , F{})  ,$
$ S_{[\omega_0 {}] } \big)$.
Then, the probability
that
$D_x^{1-\alpha} \ni \pi(\omega_0)$
is more than or equal to $1- \alpha$.
\end{enumerate}
From the ($C_2$), we assert "the statistical hypothesis test" as follows:
\rm
\begin{itemize}
\item[(D$_2$)]
[The statistical hypothesis test].
\sl
Assume that
a state $\omega_0$
satisfies that
$\pi(\omega_0)
\in
H_N
( \subseteq \Theta )
$,
where $H_N$
is called a 
"null hypothesis".
Put
\begin{align}
&
{\widehat R}_{{H_N}}^{\alpha; \Theta}
=
\bigcap_{\omega \in  \Omega \mbox{ \footnotesize such that }
\pi(\omega) \in {H_N}}
\{
E({x})
(\in
\Theta)
:
d^x_\Theta ({}E(x),
\pi(\omega )
)
\ge
\eta^\alpha_{\omega }
\}.
\label{eq13}
\intertext{and also}
&
{\widehat R}_{{H_N}}^{\alpha; X}
=
E^{-1}(
{\widehat R}_{{H_N}}^{\alpha; \Theta})
=
=
\bigcap_{\omega \in  \Omega \mbox{ \footnotesize such that }
\pi(\omega) \in {H_N}}
\{
x
(\in
X)
:
d^x_\Theta ({}E(x),
\pi(\omega )
)
\ge
\eta^\alpha_{\omega }
\}.
\label{eq114}
\end{align}
which is respectively called
\it
the $({}\alpha{})$-rejection region
of
the null hypothesis
${H_N}$.
\sl
Then,
the probability,
that
the measured value $x (\in X)$
obtained
by the measurement
${\mathsf M}_{C_0(\Omega)} \big({}{\mathsf O}:= ({}X, {\cal F} , F{}),$
%
$S_{[\omega_0]} \big)$
$($
where it should be noted that $\pi(\omega_0) \in H_N )$
satisfies the following
condition (\ref{eq15}),
is less than or equal to
$\alpha$
({}e.g., $\alpha= 0.05${}).
\begin{align}
"E(x) \in
{\widehat R}_{{H_N}}^{\alpha; \Theta}"
\mbox{  or equivalently  }
"
x
\in
{\widehat R}_{{H_N}}^{\alpha; X}"
\label{eq15}
\end{align}
\end{itemize}

\par
\vskip1.0cm
\par

%
%
%
%
%

\par
\section{ANOVA in the quantum linguistic formulation of statistics}
\rm
The arguments in this section are continued from Example 1.

\subsection{The simplest example; Student's $t$-distribution}
\par
\noindent
\bf
Example 2
\rm
[Student's t-distribution ({\it cf.} \cite{Ishi9}).
Consider the measurement
${\mathsf M}_{C_0({\mathbb R} \times {\mathbb R}_+)}$
$({\mathsf O}_G^n = ({\mathbb R}^n, {\mathcal B}_{\mathbb R}^n, {{{G}}^n}) ,$
$S_{[(\mu, \sigma)]})$
in $C_0({\mathbb R} \times {\mathbb R}_+)$
in Example 1.
Thus,
we consider that
$\Omega = {\mathbb R} \times {\mathbb R}_+$,
$X={\mathbb R}^n$.
Put
$\Theta={\mathbb R}$.
Also, define the estimator
$E:X(={\mathbb R}^n) \to \Theta(={\mathbb R})$
such that
\begin{align}
E(x)=E(x_1, x_2, \ldots , x_n )
=
\overline{\mu}(x)
=
\frac{x_1 + x_2 + \cdots + x_n}{n}
\label{eq16}
\end{align}
The quantity $\pi:\Omega(={\mathbb R} \times {\mathbb R}_+)
\to
\Theta(={\mathbb R})$
is defined by
\begin{align}
\Omega(={\mathbb R} \times {\mathbb R}_+)
\ni \omega
=
(\mu, \sigma )
\mapsto \pi (\mu, \sigma )
=
\mu
\in
\Theta(={\mathbb R})
\label{eq17}
\end{align}
\rm
Also, assume that
the $\Theta (={\mathbb R})$ has the semi-distance
$d_\Theta^x (\forall x \in X)$
such that
\begin{align}
d_\Theta^x (\theta^{(1)}, \theta^{(2)})
=
\frac{|\theta^{(1)}-\theta^{(2)}|}{\sqrt{n}{\overline{\sigma}(x)}}
=
\frac{|\theta^{(1)}-\theta^{(2)}|}{\sqrt{\overline{SS}(x)}}
\quad
\qquad
(\forall x \in X={\mathbb R}^n,
\forall \theta^{(1)}, \theta^{(2)} \in \Theta={\mathbb R}
)
\label{eq18}
\end{align}
where ${\overline{\sigma}(x)}$ is motivated by Fisher's maximum likelihood method
(see the formulas (\ref{eq4}) and (\ref{eq5})).
\par
\noindent
Define the null hypothesis $H_N$
$(\subseteq
\Theta=
{\mathbb  R}  )
)$
such that
\begin{align}
H_N=  \{\mu_0\}
\label{eq19}
\end{align}
Thus, for any
$ \omega=(\mu_0, \sigma )  ({}\in \Omega=
{\mathbb  R} \times {\mathbb R}_+ )$,
we see that
\begin{align}
&
[G^n(\{ x \in X \;:\; 
d^x_\Theta ( E(x) , \pi( \omega ) )
\ge \eta
\}
)](\omega )
\nonumber
\\
=&
[G^n(\{ x \in X \;:\; 
\frac{
|\overline{\mu}(x)- \mu_0 |}{
{{{\sqrt{\overline{SS}(x)}}}}
}
\ge \eta
\}
)](\omega )
\nonumber
\\
=
&
\frac{1}{({{\sqrt{2 \pi }\sigma{}}})^n}
\underset{
\eta \sqrt{n-1}
\le
\frac{
|\overline{\mu}(x)- \mu_0 |}{
{\sqrt{\overline{SS}(x)}}
/\sqrt{n-1}
}
}{\int \cdots \int}
\exp[{}- \frac{\sum_{k=1}^n ({}{}{x_k} - {}{\mu_0}  {})^2 
}
{2 \sigma^2}    {}] d {}{x_1} d {}{x_2}\cdots dx_n
\nonumber
\\
=
&
\frac{1}{({{\sqrt{2 \pi }{}}})^n}
\underset{
\eta^2 n({n-1})
\le
\frac{
n(\overline{\mu}(x))^2
}{
{\overline{SS}(x)}/({n-1})
}
}
{\int \cdots \int}
\exp[{}- \frac{\sum_{k=1}^n ({}{}{x_k}  {}  {})^2 
}
{2 }    {}] d {}{x_1} d {}{x_2}\cdots dx_n
\label{eq20}
\intertext{
\begin{itemize}
\item[(E$_1$)]
using the formula of Gauss integrals derived in Kolmogorov's probability theory
(also, see Remark 2 below), we finally get as follows.
\end{itemize}
}
=
&
\int^{\infty}_{
\eta^2 n({n-1})
} p_{(1,{{n}}-1) }^F(t) dt =
\alpha \;\; (\mbox{ e.g., $\alpha=0.05$})
\label{eq21}
\end{align}
where
$p_{(1,{{n}}-1) }^F$
is
the probability density function
of 
the $F$-distribution with $(1,{{n}}-1) $ degrees of freedom.
Recall the probability density function
$p_{(n_1,n_2)}^F(x)$
of
the $F$-distribution with $(n_1,n_2) $ degrees of freedom
is represented as follows:
\begin{align}
p_{(n_1,n_2)}^F(t)
=
\frac{1}{B(n_1/2, n_2/2)}
\Big(\frac{n_1}{n_2} \Big)^{n_1/2}
\frac{t^{(n_1-2)/2}}{(1+n_1t/n_2)^{(n_1+n_2)/2}}
\qquad
(t \ge 0)
\label{eq22}
\end{align}
where $B(\cdot, \cdot)$ is the Beta function.
Define the $\alpha$-point $F_{n_1, \alpha}^{n_2}$
$( > 0)$
such that
\begin{align}
\int^{\infty}_{F_{n_1, \alpha}^{n_2} } p_{(n_1,n_2) }^F (t) dt =\alpha
\qquad
(0 < \alpha \ll 1. \mbox{  e.g., } \alpha=0.05)
\label{eq23}
\end{align}
Thus, it suffices to put
\begin{align}
{\eta^2 n({{n}}-1) }{ }
={F_{n-1, \alpha}^{1} } 
\label{eq24}
\end{align}
And thus,
\begin{align}
(\eta^\alpha_{\omega})^2 =
\frac{{F_{n-1, \alpha}^{1} }}{n(n-1)}
\label{eq25}
\end{align}
\par
\noindent

Therefore,
we get ${\widehat R}_{{H_N}}^{\alpha; \Theta}$(
or
${\widehat R}_{H_N}^{\alpha; X}$;
\it
the $({}\alpha{})$-rejection region
of
$H_N(=\{\mu_0\})$
\rm
)
as follows:
\begin{align}
{\widehat R}_{{H_N}}^{\alpha; \Theta}
&
=
\bigcap_{\omega =(\mu, \sigma ) \in  \Omega (={\mathbb R} \times {\mathbb R}_+) \mbox{ \footnotesize such that }
\pi(\omega)= \mu \in {H_N}(=\{\mu_0\})}
\{
E(x)
(\in
\Theta)
:
\;\;
d^x_\Theta ({}E(x),
\pi(\omega )
)
\ge
\eta^\alpha_{\omega }
\}
\nonumber
\\
&
=
\{\overline{\mu}(x) \in \Theta(={\mathbb R})
\;:\;
\frac{
|\overline{\mu}(x)- \mu_0 |}{
\sqrt{{\overline{SS}(x)}}
}
\ge
\eta_\omega^\alpha
\}
=
\{\overline{\mu}(x) \in \Theta(={\mathbb R})
\;:\;
\frac{
|\overline{\mu}(x)- \mu_0 |}{
\overline{\sigma}(x)
}
\ge
\eta_\omega^\alpha \sqrt{n}
\}
\nonumber
\\
&
=
\{\overline{\mu}(x) \in \Theta(={\mathbb R})
\;:\;
\frac{
|\overline{\mu}(x)- \mu_0 |}{
\overline{\sigma}(x)
}
\ge
\sqrt{\frac{F_{n-1, \alpha}^1}{n-1}}
\;\;\}
\nonumber
\\
&
=
\{\overline{\mu}(x) \in \Theta(={\mathbb R})
\;:\;
\mu_0 \le 
\overline{\mu}(x)
-
{{\overline{\sigma}(x)}}
\sqrt{\frac{F_{n-1, \alpha}^1}{n-1}}
\mbox{ or }
\overline{\mu}(x)
+
{{\overline{\sigma}(x)}}
\sqrt{\frac{F_{n-1, \alpha}^1}{n-1}}
\le \mu_0
\}
\label{eq26} 
\end{align}
and
\begin{align}
{\widehat R}_{H_N}^{\alpha; X}
&=
E^{-1}({\widehat R}_{{H_N}}^{\alpha; \Theta})
\nonumber
\\
&
=
\{x \in X(={\mathbb R}^n)
\;:\;
\mu_0 \le 
\overline{\mu}(x)
-
{{\overline{\sigma}(x)}}
\sqrt{\frac{F_{n-1, \alpha}^1}{n-1}}
\mbox{ or }
\overline{\mu}(x)
+
{{\overline{\sigma}(x)}}
\sqrt{\frac{F_{n-1, \alpha}^1}{n-1}}
\le \mu_0
\}
\label{eq27} 
\end{align}
Therefore, the statistical hypothesis test (D$_2$) in Theorem 1 is applicable.

\rm
\vskip0.5cm
\par
\noindent
\bf
Remark 2
\rm
[Kolmogorov's probability theory ({\it cf}. \cite{Kolm}).
There are several derivations of (\ref{eq21}) from (\ref{eq20}).
Of course,
the formula (\ref{eq21}) can be directly derived from the (\ref{eq20}),
though the calculation is not easy.
However, Kolmogorov's probability theory is useful for the derivation of (((1).
Here, let us remark it as follows.
Consider the probability space
$({\mathbb R}^n, {\mathcal B}_{{\mathbb R}^n},P )$,
where
$P(\Xi) =[G^n(\Xi)](\mu, \sigma)$
$(\forall \Xi \in
{\mathcal B}_{{\mathbb R}^n})$. And, for each $k=1,2, \ldots, n$,
consider a random variable $X_k : {\mathbb R}^n \to {\mathbb R}$ such that
$$
X_k (x) =X_k(x_1,x_2, \ldots, x_n) = x_k
\quad
(\forall x=(x_1,x_2, \ldots, x_n) \in {\mathbb R}^n
$$
It is clear that
random variables $\{X_k \}_{k=1,2,\ldots, n}$ are independent
with the normal distribution $N(\mu, \sigma^2)$.
Define the random variables
$\overline{\mu}:{\mathbb R}^n \to {\mathbb R}$
and
$\overline{\sigma}:{\mathbb R}^n \to {\mathbb R}$,
which are also independent
({\it cf}. Cochran's theorem, etc.).
And thus, we can easily show that
the random variable
$\frac{
n{\overline{\mu}(x) }^2
}{
{{\overline{\sigma}(x)}^2}}$
has the $F$-distribution 
with $(1,{{n}}-1) $ degrees of freedom. 
Also,
the random variable
$\frac{
{\overline{\mu}(x)- \mu }
}{
{{\overline{\sigma}(x)}}/\sqrt{n-1}}$
has the student's $t$-distribution.
Therefore,
Kolmogorov's probability theory
provides a useful calculation method
to quantum language.
We never consider that
Kolmogorov's probability theory gives a foundation to statistics.
However,
it is certain that mathematical theories
(particularly,
the theory of probability
and the theory of operator algebra
({\it cf}. \cite{Saka})
)
are
indispensable for quantum language.
For completeness,
again note that
Kolmogorov's probability theory
is merely used in order to calculate
multi-dimensional Gauss integrals
throughout this paper
({\it cf}.
the items
(E$_1$)-(E$_4$)
in Examples 2-5).


\subsection{The one-way ANOVA}
\par
\noindent
\bf
Example 3
\rm
[The one-way ANOVA].
For each $i=1,2, \cdots , a$,
a natural number $n_i$ is determined. And put
${{n}}=\sum_{i=1}^a n_i$.
As one of generalizations of Example 2,
we consider a kind of normal observable
${\mathsf O}_G^{{{n}}} = (X(\equiv {\mathbb R}^{{{n}}}), {\mathcal B}_{\mathbb R}^{{{n}}}, {{{G}}^{{{n}}}} )$ 
in $C_0(\Omega ( \equiv ({\mathbb R}^a \times {\mathbb R}_+))$
as follows:
\par
\noindent
\begin{align}
&
[{{{G}}}^{{{n}}}
(
\widehat{\Xi})
]
({}\omega{})
=
\frac{1}{({{\sqrt{2 \pi }\sigma{}}})^{{{n}}}}
\underset{\widehat{\Xi}
}{\int \cdots \int}
\exp[{}- \frac{\sum_{i=1}^a \sum_{k=1}^{n_i}  ({}{}{x_{ik}} - {}{\mu}_i  {})^2 
}
{2 \sigma^2}   {}] 
\bigtimes_{i=1}^a \bigtimes_{k=1}^{n_i}
d {}{x_{ik}} 
\label{eq28}
\\
&
\qquad
(
\forall \omega =(\mu_1, \mu_2, \ldots, \mu_a, \sigma)
 \in \Omega = {\mathbb R}^a \times {\mathbb R}_+ ,
 \widehat{\Xi} \in {\mathcal B}_{\mathbb R}^{{{n}}})
\nonumber
\end{align}
Put
\begin{align}
\alpha_i= \mu_i - \frac{\sum_{i=1}^a \mu_i }{a}
\qquad
(\forall i=1,2, \ldots, a )
\label{eq29} 
\end{align}
and
\begin{align}
\Theta = {\mathbb R}^a
\label{eq30} 
\end{align}
and define the map $\pi : \Omega \to \Theta $
such that
\begin{align}
\Omega = {\mathbb R}^a \times {\mathbb R}_+
\ni
\omega =(\mu_1, \mu_2, \ldots, \mu_a, \sigma)
\mapsto
\pi(\omega)
=
(\alpha_1, \alpha_2, \ldots, \alpha_a)
\in
\Theta = {\mathbb R}^a
\label{eq31} 
\end{align}
Define the null hypothesis
$H_N ( \subseteq \Theta = {\mathbb R}^a)$
such that
\begin{align}
H_N
&
=
\{
(\alpha_1, \alpha_2, \ldots, \alpha_a)
\in
\Theta = {\mathbb R}^a
\;:\;
\alpha_1=\alpha_2= \ldots= \alpha_a= \alpha
\}
\nonumber
\\
&
=
\{ ( \overbrace{0, 0, \ldots, 0}^{a} ) \}
\label{eq32} 
\end{align}
since it clearly holds that
$
"\mu_1=\mu_2=\ldots=\mu_a"
\Leftrightarrow
"\alpha_1=\alpha_2=\ldots=\alpha_a=0"
$.
\par
\noindent
Put
\begin{align}
&
\|
\theta^{(1)}- \theta^{(2)}
\|_\Theta
=
\sqrt{
\sum_{i=1}^a n_i \Big(\theta_{i}^{(1)} - \theta_{i}^{(2)} 
\Big)^2
}
\label{eq33} 
\\
&
\qquad
(\forall \theta^{(\ell)} =( \theta_1^{(\ell)}, \theta_2^{(\ell)}, \ldots, \theta_a^{(\ell)} )
\in {\mathbb R}^{a},
\;
\ell=1,2
)
\nonumber
\end{align}
which is the weighted Euclidean norm in ${\mathbb R}^{n_i}$.
\par
\noindent
Also, put
\begin{align}
&X={\mathbb R}^{{{n}}}
\ni
x
=
((x_{ik})_{ k=1,2, \ldots, n_i})_{i=1,2,\ldots,a}
\nonumber
\\
& x_{i \cdot} =\frac{\sum_{k=1}^{n_i} x_{ik}}{n_i}, \qquad
x_{ \cdot \cdot} =\frac{\sum_{i=1}^a \sum_{k=1}^{n_i}x_{ik}}{{{n_i}}}, \quad
\label{eq34} 
\end{align}
\par
\noindent
According to Fisher's maximum likelihood method,
define and calculate $\overline{\sigma}(x) (=
\sqrt{
\frac{{\overline{SS}}(x)}{n}
}
)$
concerning (\ref{eq28})
as follows. For each $x \in X={\mathbb R}^{{{n}}}$,
\begin{align}
&
{\overline{SS}}(x)
=
{\overline{SS}}(((x_{ik})_{\; k=1,2, \ldots, {n_i} })_{i=1,2, \ldots, a\;}
)
\nonumber
\\
=
&
\sum_{i=1}^a \sum_{k=1}^{n_i} (x_{ik} - x_{i \cdot})^2
\nonumber
\\
=
&
\sum_{i=1}^a \sum_{k=1}^{n_i} (x_{ik} - \frac{\sum_{k=1}^{n_i} x_{i k}}{n_i})^2
\nonumber
\\
=
&
\sum_{i=1}^a \sum_{k=1}^{n_i}  ((x_{ik}-\mu_i) - \frac{\sum_{k=1}^{n_i} (
x_{i k}-\mu_i)}{n_i})^2
\qquad
\nonumber
\\
=
&
{\overline{SS}}(((x_{ik}- \mu_{i})_{\; k=1,2, \ldots, {n_i} })_{i=1,2, \ldots, a\;}
)
\label{eq35}
\end{align}
\par
\noindent
And, for each $x \in X = {\mathbb R}^{{{n}}}$, define
the semi-distance
$d_\Theta^x$ in $\Theta$ such that
\begin{align}
&
d_\Theta^x (\theta^{(1)}, \theta^{(2)})
=
\frac{\|\theta^{(1)}- \theta^{(2)} \|_\Theta}{
\sqrt{{\overline{SS}}(x)
}
}
\qquad
(\forall \theta^{(1)}, \theta^{(2)} \in \Theta )
 ).
\label{eq36}
\end{align}

\par
\noindent
Further define the estimator
$E: X(={\mathbb R}^{{{n}}}) \to \Theta(={\mathbb R}^{a} )$ such that
\begin{align}
E(x)
=
&
E(
(x_{ik})_{i=1,2,\ldots,a, k=1,2, \ldots, n}
)
\nonumber
\\
=
&
\Big(
\frac{\sum_{k=1}^{n_i} x_{1k}}{n}
-
\frac{ \sum_{i=1}^a \sum_{k=1}^{n_i}   x_{ik}}{{{n}}}
,
\frac{\sum_{k=1}^{n_i} x_{2k}}{n}
-
\frac{ \sum_{i=1}^a \sum_{k=1}^{n_i}   x_{ik}}{{{n}}},
\ldots,
\frac{\sum_{k=1}^{n_i} x_{ak}}{n}
-
\frac{ \sum_{i=1}^a \sum_{k=1}^{n_i}   x_{ik}}{{{n}}}
\Big)
\nonumber
\\
=
&
\Big(
\frac{\sum_{k=1}^{n_i} x_{ik}}{n}
-
\frac{ \sum_{i=1}^a \sum_{k=1}^{n_i}   x_{ik}}{{{n}}}
\Big)_{i=1,2, \ldots, a }
=
(x_{i \cdot} - x_{\cdot \cdot })_{i=1,2, \ldots, a }
\label{eq37}
\end{align}
Hence, we see that
\begin{align}
&
\| E(x) - \pi (\omega )\|^2_\Theta
\nonumber
\\
=
&
||
\Big(
\frac{\sum_{k=1}^{n_i} x_{ik}}{n}
-
\frac{ \sum_{i=1}^a \sum_{k=1}^{n_i}   x_{ik}}{{{n}}}
\Big)_{i=1,2, \ldots, a }
-
(
\alpha_i
)_{i=1,2, \ldots, a }
||_\Theta^2
\nonumber
\\
=
&
||
\Big(
\frac{\sum_{k=1}^{n_i} x_{ik}}{n}
-
\frac{ \sum_{i=1}^a \sum_{k=1}^{n_i}   x_{ik}}{{{n}}}
-
(\mu_i - \frac{\sum_{i=1}^a \mu_i }{a})
\Big)_{i=1,2, \ldots, a }
||_\Theta^2
\label{eq38}
\intertext{and thus, if the null hypothesis $H_N$ is assumed
(i.e.,
$\mu_i-\frac{\sum_{k=1}^a\mu_i}{a}=\alpha_i =0 (i=1,2,\ldots, a )$), we see
}
=
&
||
\Big(
\frac{\sum_{k=1}^{n_i} x_{ik}}{n}
-
\frac{ \sum_{i=1}^a \sum_{k=1}^{n_i}   x_{ik}}{{{n}}}
\Big)_{i=1,2, \ldots, a }
||_\Theta^2
=
\sum_{i=1}^a n_i (x_{i \cdot} - x_{\cdot \cdot})^2
\label{eq39}
\end{align}
Thus,
for any
$ \omega=((\mu_{ik})_{i=12,\ldots,a, \;k=1,2, \ldots, n }, \sigma )  ({}\in \Omega=
{\mathbb  R}^{{{n}}} \times {\mathbb R}_+ )$,
define
the positive number
$\eta^\alpha_{\omega}$
$({}> 0)$
such that:
\begin{align}
\eta^\alpha_{\omega}
=
\inf
\{
\eta > 0:
[G^{{{n}}}({}E^{-1} ({}
{{\rm Ball}^C_{d_\Theta^{x}}}(\pi(\omega) ; \eta{}))](\omega )
\ge \alpha
\}
\label{eq40}
\end{align}
where
\begin{align}
{\rm Ball}^C_{d_\Theta^{x}}(\pi(\omega) ; \eta{})
=\{ \theta \in \Theta
\;:\;
d_\Theta^{x} ( \pi(\omega ) , \theta ) > \eta \}
\label{eq41}
\end{align}
Assume the null hypothesis $H_N$
(i.e.,
$\mu_i-\frac{\sum_{k=1}^a\mu_i}{a}=\alpha_i =0 (i=1,2,\ldots, a )$).
Now let us calculate the $\eta^\alpha_{\omega}$ as follows:
\begin{align}
&
E^{-1}({{\rm Ball}^C_{d_\Theta^{x} }}(\pi(\omega) ; \eta ))
=\{ x \in X = {\mathbb R}^{{{n}}}
\;:\;
d_\Theta^x (E(x), \pi(\omega ))
>
\eta
\}
\nonumber
\\
=
&
\{ x \in X = {\mathbb R}^{{{n}}}
\;:\;
\frac{
\| E(x)- \pi(\omega) \|^2_\Theta
}{{{\overline{SS}}(x) }}
=
\frac{
\sum_{i=1}^a n_i (
x_{i \cdot}
- x_{\cdot \cdot}
)^2}{
\sum_{i=1}^a \sum_{k=1}^{n_i} (x_{ik} - x_{i \cdot})^2
}
>
\eta^2
\}
\label{eq42}
\end{align}


\par
\noindent
That is,
for any $\omega
=(\mu_1, \mu_2, \ldots, \mu_a, \sigma) \in \Omega={\mathbb R}^{a} \times {\mathbb R}_+$
such that
$\pi( \omega ) (=
(\alpha_1, \alpha_2, \ldots, \alpha_a)
)\in H_N (=\{0,0, \ldots, 0)\})$,

\par
\noindent
\begin{align}
&
[{{{G}}}^{{{n}}}
(
E^{-1}({{\rm Ball}^C_{d_\Theta^{x} }}(\pi(\omega) ; \eta ))
)
({}\omega{})
\nonumber
\\
=
&
\frac{1}{({{\sqrt{2 \pi }\sigma{}}})^{{{n}}}}
\underset{
\frac{
\sum_{i=1}^a n_i (
x_{i \cdot}
- x_{\cdot \cdot}
)^2}{
\sum_{i=1}^a \sum_{k=1}^{n_i} (x_{ik} - x_{i \cdot})^2
}
> \eta^2
}{\int \cdots \int}
\exp[{}- \frac{ \sum_{i=1}^a \sum_{k=1}^{n_i}   ({}{}{x_{ik}} - {}{\mu_i}  {})^2 
}
{2 \sigma^2}   {}] 
\bigtimes_{i=1}^a \bigtimes_{k=1}^{n_i}
d {}{x_{ik}} 
\nonumber
\\
=
&
\frac{1}{({{\sqrt{2 \pi }{}}})^{{{n}}}}
\underset{
\frac{
(\sum_{i=1}^a n_i(
x_{i \cdot}
- x_{\cdot \cdot}
)^2 /(a-1)}{
(\sum_{i=1}^a \sum_{k=1}^{n_i} (x_{ik} - x_{i \cdot})^2)/({{n}}-a)
}
> \frac{\eta^2 ({{n}}-a) }{ (a-1)}
}
{\int \cdots \int}
\exp[{}- \frac{ \sum_{i=1}^a \sum_{k=1}^{n_i}   ({}{}{x_{ik}}  {})^2 
}
{2 }   {}] 
\bigtimes_{i=1}^a \bigtimes_{k=1}^{n_i}
d {}{x_{ik}} 
\label{eq43}
\intertext{
\begin{itemize}
\item[(E$_2$)]
using the formula of Gauss integrals derived in Kolmogorov's probability theory
(also, recall Remark 2), we finally get as follows.
\end{itemize}
}
=
&
\int^{\infty}_{
\frac{\eta^2 ({{n}}-a) }{ (a-1)}
} p_{(a-1,{{n}}-a) }^F(t) dt =
\alpha \;\; (\mbox{ e.g., $\alpha$=0.05})
\label{eq44}
\end{align}
where
$p_{(a-1,{{n}}-a) }^F$
is
the probability density function 
the $F$-distribution with $(a-1,{{n}}-a) $ degrees of freedom.
Thus, it suffices to put
\begin{align}
\frac{\eta^2 ({{n}}-a) }{ (a-1)}
={F_{n-a, \alpha}^{a-1} } 
(=\mbox{"$\alpha$-point"})
\label{eq45}
\end{align}
And thus we see,
\begin{align}
(\eta^\alpha_{\omega})^2 =
{F_{n-a, \alpha}^{a-1} }
(a-1)/(n-a)
\label{eq46}
\end{align}
\par
\noindent
Therefore,
we get ${\widehat R}_{\widehat{x}}^{\alpha; \Theta}$
(or,
${\widehat R}_{\widehat{x}}^{\alpha; X}$;
the $({}\alpha{})$-rejection region
of
$H_N =\{(0.0. \ldots, 0)\}( \subseteq \Theta= {\mathbb R}^a)$
)
as follows:
\begin{align}
{\widehat R}_{{H_N}}^{\alpha; \Theta}
&
=
\bigcap_{\omega =((\mu_i)_{i=1}^a, \sigma ) \in  \Omega (={\mathbb R}^a \times {\mathbb R}_+ ) \mbox{ \footnotesize such that }
\pi(\omega)= (\mu)_{i=1}^a \in {H_N}=\{(0,0,\ldots,0)\}}
\{
E({x})
(\in
\Theta)
:
d_\Theta^{x}  ({}E({x}),
\pi(\omega))
\ge
\eta^\alpha_{\omega }
\}
\nonumber
\\
&
=
\{
E({x})
(\in
\Theta)
:
\frac{
(\sum_{i=1}^a n_i (
x_{i \cdot}
- x_{\cdot \cdot}
)^2) /(a-1)}{
(\sum_{i=1}^a  \sum_{k=1}^{a_i} (x_{ik} - x_{i \cdot})^2))/({{n}}-a)
}
\ge
{F_{n-a, \alpha}^{a-1} }
\}
\label{eq47}
\end{align}
Thus,
\begin{align}
{\widehat R}_{\widehat{x}}^{\alpha; X}
=
E^{-1}({\widehat R}_{H_N}^{\alpha;\Theta})
=
\{
 x \in X \;:\;
{
\frac{
(\sum_{i=1}^a n_i (
x_{i \cdot}
- x_{\cdot \cdot}
)^2 )/(a-1)}{
(\sum_{i=1}^a \sum_{k=1}^{n_i} (x_{ik} - x_{i \cdot})^2)/({{n}}-a)
}
\ge
{F_{n-a, \alpha}^{a-1} }
}
\}
\label{eq48}
\end{align}
Therefore, the statistical hypothesis test (D$_2$) in Theorem 1 is applicable.

\rm
\vskip0.5cm
\par
\par
\noindent

\par
\subsection{The two-way ANOVA}
As one of generalizations of Example 2,
we consider a kind of observable
${\mathsf O}_G^{abn} = (X(\equiv {\mathbb R}^{abn}), {\mathcal B}_{\mathbb R}^{abn}, {{{G}}^{abn}} )$ 
in $C_0(\Omega ( \equiv ({\mathbb R}^{ab} \times {\mathbb R}_+))$.

\par
\noindent
\rm
\par
\par
\noindent
Put
\begin{align}
&
[{{{G}}}^{abn}
(\widehat{\Xi})]
({}\omega{})
\nonumber
\\
=
&
\frac{1}
{({
{\sqrt{2 \pi }
\sigma}
})^{abn}}
\underset{
\widehat{\Xi}
}
{\int \cdots \int}
\exp[- 
\frac{
\sum_{i=1}^a \sum_{j=1}^b \sum_{k=1}^n  (x_{ijk} - 
\mu_{ij}
)^2 
}{2 \sigma^2}
] 
\bigtimes_{k=1}^n
\bigtimes_{j=1}^b
\bigtimes_{i=1}^a
d{x_{ijk} }
\label{eq49}
\\
&
\qquad
(
\forall \omega =((\mu_{ij})_{i=1,2, \ldots, a,j=1,2, \ldots, b}, \sigma)
 \in \Omega = {\mathbb R}^{ab+1} \times {\mathbb R}_+ ,
 \widehat{\Xi} \in {\mathcal B}_{\mathbb R}^{abn}
 )
\nonumber
\end{align}
\par
\noindent
Put
\begin{align}
\mu_{ij}
&
=
\overline{\mu} (=
\mu_{\cdot \cdot }
= \frac{\sum_{i=1}^a \sum_{j=1}^b \mu_{ij} }{ab})
\nonumber
\\
&
\quad
+
\alpha_i
(= \mu_{i \cdot} - \mu_{\cdot \cdot }
=
\frac{\sum_{j=1}^b  \mu_{ij} }{b}
-
\frac{\sum_{i=1}^a \sum_{j=1}^b \mu_{ij} }{ab}
)
\nonumber
\\
&
\quad
+
\beta_j
(=\mu_{\cdot j } -  \mu_{\cdot \cdot }
=
\frac{\sum_{i=1}^a  \mu_{ij} }{a}
-
\frac{\sum_{i=1}^a \sum_{j=1}^b \mu_{ij} }{ab}
)
\nonumber
\\
&
\quad
+
{(\alpha \beta)}_{ij}
(=\mu_{ij} -\mu_{i \cdot}-\mu_{ \cdot j}+ \mu_{\cdot \cdot } )
\label{eq50}
\end{align}
Put
\begin{align}
&X={\mathbb R}^{abn}
\ni
x
=
(x_{ijk})_{i=1,2,\ldots,a,\; j=1,2,\ldots,b, \; k=1,2, \ldots, n}
\nonumber
\\
&
x_{ij\cdot}= \frac{\sum_{k=1}^{n}x_{ijk}}{n}, \quad
x_{i \cdot \cdot} =\frac{\sum_{j=1}^b \sum_{k=1}^n x_{ijk}}{bn}, \quad
x_{ \cdot j \cdot} =\frac{\sum_{i=1}^a \sum_{k=1}^n x_{ijk}}{an}, \quad
\nonumber
\\
&
x_{ \cdot \cdot \cdot} =\frac{\sum_{i=1}^a \sum_{j=1}^b \sum_{k=1}^n x_{ijk}}{abn} \quad
\label{eq51}
\end{align}

\par
\noindent
\bf
Example 4
\rm
[The null hypothesis such that $\alpha_1=\alpha_2=\cdots=\alpha_a=0$].
\rm
Define the $\pi : \Omega \to \Theta $
such that
\begin{align}
\Omega = {\mathbb R}^{ab+1} \times {\mathbb R}_+
\ni
\omega =((\mu_{ij})_{i=1,2, \ldots, a,j=1,2, \ldots, b}, \sigma)
\mapsto
\pi_1(\omega)
=
(\alpha_i)_{i=1}^a
\in
\Theta = {\mathbb R}^a
\label{eq52}
\end{align}

%
Put
\begin{align}
\Theta = {\mathbb R}^a
\label{eq53}
\end{align}
and define the $\pi : \Omega \to \Theta $
such that
\begin{align}
\Omega = {\mathbb R}^a \times {\mathbb R}_+
\ni
\omega =(\mu_1, \mu_2, \ldots, \mu_a, \sigma)
\mapsto
\pi(\omega)
=
(\alpha_1, \alpha_2, \ldots, \alpha_a)
\in
\Theta = {\mathbb R}^a
\label{eq54}
\end{align}
Define the null hypothesis
$H_N ( \subseteq \Theta = {\mathbb R}^a)$
such that
\begin{align}
H_N
&
=
\{
(\alpha_1, \alpha_2, \ldots, \alpha_a)
\in
\Theta = {\mathbb R}^a
\;:\;
\alpha_1=\alpha_2= \ldots= \alpha_a= \alpha
\}
\nonumber
\\
&
=
\{ ( \overbrace{0, 0, \ldots, 0}^{a} ) \}
\label{eq55}
\end{align}
That is because
\begin{align}
a \alpha=\sum_{i=1}^a \alpha_i
=\sum_{i=1}^a 
(\mu_{i \cdot} - \mu_{\cdot \cdot })
=
\frac{\sum_{i=1}^a\sum_{j=1}^b  \mu_{ij} }{b}
-
\sum_{i=1}^a \frac{\sum_{i=1}^a \sum_{j=1}^b \mu_{ij} }{ab}
=0
\label{eq56}
\end{align}
Put
\begin{align}
&
\|
\theta^{(1)}- \theta^{(2)}
\|_\Theta
=
\sqrt{
\sum_{i=1}^a \Big(\theta_{i}^{(1)} - \theta_{i}^{(2)} 
\Big)^2
}
\label{eq57}
\\
&
\qquad
(\forall \theta^{(\ell)} =( \theta_1^{(i)}, \theta_2^{(\ell)}, \ldots, \theta_a^{(\ell)} )
\in {\mathbb R}^{a},
\;
\ell=1,2
)
\nonumber
\end{align}


\par
\noindent
Motivated by Fisher's maximum likelihood method,
define and calculate $\overline{\sigma}(x) \Big(=
\sqrt{{\overline{SS}}(x)/(abn)}
\Big)$
as follows.
\begin{align}
&
{\overline{SS}}(x)=
{\overline{SS}}((x_{ijk})_{i=1,2, \ldots, a,\;\; j=1,2, \ldots, b,k=1,2, \ldots, n })
\nonumber
\\
:=
&
{
\sum_{i=1}^a \sum_{j=1}^b \sum_{k=1}^n 
(x_{ijk} - x_{ij \cdot})^2
}
=
{
\sum_{i=1}^a \sum_{j=1}^b \sum_{k=1}^n (x_{ijk} - \frac{\sum_{k=1}^n x_{ij k}}{n})^2
}
\nonumber
\\
=
&
{
\sum_{i=1}^a \sum_{j=1}^b \sum_{k=1}^n ((x_{ijk}-\mu_{ij}) - \frac{\sum_{k=1}^n (
x_{i jk}-\mu_{ij})}{n})^2
}
\nonumber
\\
=
&
{\overline{SS}}(((x_{ijk}- \mu_{ij})_{i=1,2, \ldots, a,\;\; j=1,2, \ldots, b })_{k=1,2, \cdots, n})
\label{eq58}
\end{align}

\par
\noindent
Define
the semi-distance
$d_\Theta^x$ in
$\Theta = {\mathbb R}^a$ such that
\begin{align}
&
d_\Theta^x (\theta^{(1)}, \theta^{(2)})
=
\frac{\|\theta^{(1)}- \theta^{(2)} \|_\Theta}{
\sqrt{{\overline{SS}}(x)}
}
\qquad
(\forall \theta^{(1)}, \theta^{(2)}
\in \Theta={\mathbb R}^{a},
 \forall x \in X={\mathbb R}^{abn} )
\label{eq59}
\end{align}

\par
\noindent
Define the estimator
$E: X(={\mathbb R}^{abn}) \to \Theta(={\mathbb R}^{a} )$ such that
\begin{align}
E(x)
=
\Big(\frac{\sum_{j=1}^b\sum_{k=1}^n x_{ijk}}{bn}
-
\frac{\sum_{i=1}^a \sum_{j=1}^b \sum_{k=1}^n x_{ijk}}{abn}
\Big)_{i=1,2,\ldots, a }
=
\Big( x_{i \cdot \cdot }
-
x_{\cdot \cdot \cdot }
\Big)_{i=1,2,\ldots, a }
\label{eq60}
\end{align}

\par
\noindent
Hence
\begin{align}
&
\| E(x) - \pi (\omega )\|^2_\Theta
\nonumber
\\
=
&
||
\Big(\frac{\sum_{j=1}^b\sum_{k=1}^n x_{ijk}}{bn}
-
\frac{\sum_{i=1}^a \sum_{j=1}^b \sum_{k=1}^n x_{ijk}}{abn}
\Big)_{i=1,2,\ldots, a }
-
\Big(
\alpha_i
\Big)_{i=1,2, \ldots, a }
||_\Theta^2
\nonumber
\\
=
&
||
\Big(\frac{\sum_{j=1}^b\sum_{k=1}^n x_{ijk}}{bn}
-
\frac{\sum_{i=1}^a \sum_{j=1}^b \sum_{k=1}^n x_{ijk}}{abn}
\Big)_{i=1,2,\ldots, a }
-
\Big(
\frac{\sum_{j=1}^b  \mu_{ij} }{b}
-
\frac{\sum_{i=1}^a \sum_{j=1}^b \mu_{ij} }{ab}
\Big)_{i=1,2, \ldots, a }
||_\Theta^2
\nonumber
\\
=
&
||
\Big(
\frac{\sum_{k=1}^n \sum_{j=1}^b(x_{ijk}-\mu_{ij})}{bn}
-
\frac{\sum_{i=1}^a \sum_{j=1}^b \sum_{k=1}^n (x_{ijk}-\mu_{ij})}{abn}
\Big)_{i=1,2, \ldots, a }
||_\Theta^2
\label{eq61}
\intertext{and thus, if the null hypothesis $H_N$ is assumed
(i.e.,
$\mu_{i \cdot} -\mu_{\cdot \cdot} =\alpha_i=0$
$(\forall i=1,2,\ldots, a )$
)
}
=
&
||
\Big(
\frac{\sum_{k=1}^n \sum_{j=1}^bx_{ijk}}{bn}
-
\frac{\sum_{i=1}^a \sum_{j=1}^b \sum_{k=1}^n x_{ijk}}{abn}
\Big)_{i=1,2, \ldots, a }
||_\Theta^2
=
\sum_{i=1}^a(x_{ij \cdot} - x_{\cdot \cdot \cdot})^2
\label{eq62}
\end{align}

Thus,
for any
$ \omega=(\mu_1, \mu_2 )  ({}\in \Omega=
{\mathbb  R} \times {\mathbb R} )$,
define
the positive number
$\eta^\alpha_{\omega}$
$({}> 0)$
such that:
\begin{align}
\eta^\alpha_{\omega}
=
\inf
\{
\eta > 0:
[G({}E^{-1} ({}
{{\rm Ball}^C_{d_\Theta^{x}}}(\pi(\omega) ; \eta{}))](\omega )
\ge \alpha
\}
\label{eq63}
\end{align}

Assume the null hypothesis $H_N$.
Now let us calculate the $\eta^\alpha_{\omega}$ as follows:
\begin{align}
&
E^{-1}({{\rm Ball}^C_{d_\Theta^{x} }}(\pi(\omega) ; \eta ))
=\{ x \in X = {\mathbb R}^{abn}
\;:\;
d_\Theta^x (E(x), \pi(\omega ))
>
\eta
\}
\nonumber
\\
=
&
\{ x \in X = {\mathbb R}^{abn}
\;:\;
\frac{
abn \sum_{i=1}^a \sum_{j=1}^b(
x_{ij \cdot}
- x_{\cdot \cdot \cdot}
)^2}{
\sum_{i=1}^a \sum_{j=1}^b\sum_{k=1}^n (x_{ijk} - x_{ij \cdot})^2
}
>
\eta
\}
\label{eq64}
\end{align}

\par
\noindent
That is,
for any $\omega
=((\mu_{ij})_{i=1,2,\ldots,a, \;j=1,2,\ldots,b},\;, \sigma) \in \Omega$
such tht
$\pi( \omega ) (=
(\alpha_1, \alpha_2, \ldots, \alpha_a)
)\in H_N (=\{0,0, \ldots, 0)\})$,

\par
\noindent
\begin{align}
&
[{{{G}}}^{abn}
(
E^{-1}({{\rm Ball}^C_{d_\Theta^{x} }}(\pi(\omega) ; \eta ))
)
({}\omega{})
\nonumber
\\
=
&
\frac{1}
{({
{\sqrt{2 \pi }
\sigma}
})^{abn}}
\underset{
E^{-1}({{\rm Ball}^C_{d_\Theta^{x} }}(\pi(\omega) ; \eta ))
}
{\int \cdots \int}
\exp[- 
\frac{
\sum_{i=1}^a \sum_{j=1}^b \sum_{k=1}^n  (x_{ijk} - 
\mu_{ij}
)^2 
}{2 \sigma^2}
] 
\bigtimes_{k=1}^n
\bigtimes_{j=1}^b
\bigtimes_{i=1}^a
d{x_{ijk} }
\nonumber
\\
=
&
\frac{1}
{({
{\sqrt{2 \pi }
\sigma}
})^{abn}}
\underset{
\frac{
abn \sum_{i=1}^a \sum_{j=1}^b(
x_{ij \cdot}
- x_{\cdot \cdot \cdot}
)^2}{
\sum_{i=1}^a \sum_{j=1}^b\sum_{k=1}^n (x_{ijk} - x_{ij \cdot})^2
}
>
\eta^2}
{\int \cdots \int}
\exp[- 
\frac{
\sum_{i=1}^a \sum_{j=1}^b \sum_{k=1}^n  (x_{ijk} - 
\mu_{ij}
)^2 
}{2 \sigma^2}
] 
\bigtimes_{k=1}^n
\bigtimes_{j=1}^b
\bigtimes_{i=1}^a
d{x_{ijk} }
\nonumber
\\
=
&
\frac{1}
{({
{\sqrt{2 \pi }}
})^{abn}}
\underset{
\frac{
\frac{
\sum_{i=1}^a \sum_{j=1}^b(
x_{ij \cdot}
- x_{\cdot \cdot \cdot}
)^2)}{(a-1)}
}{
\frac{
\sum_{i=1}^a \sum_{j=1}^b\sum_{k=1}^n (x_{ijk} - x_{ij \cdot})^2
}{ab(n-1)}
}
>
\frac{\eta^2 (ab(n-1))}{abn(a-1)}
}
{\int \cdots \int}
\exp[- 
\frac{
\sum_{i=1}^a \sum_{j=1}^b \sum_{k=1}^n  (x_{ijk}
)^2 
}{2 }
] 
\bigtimes_{k=1}^n
\bigtimes_{j=1}^b
\bigtimes_{i=1}^a
d{x_{ijk} }
\label{eq65}
\intertext{
\begin{itemize}
\item[(E$_3$)]
using the formula of Gauss integrals derived in Kolmogorov's probability theory
(also, recall Remark 2), we finally get as follows.
\end{itemize}
}
=
&
\int^{\infty}_{\frac{\eta^2 (n-1)}{n (a-1)}} p_{(a-1,ab(n-1)) }^F(t) dt =
\alpha
(\mbox{e.g., } \alpha=0.05)
\label{eq66}
\end{align}
where
$p_{(a-1,ab(n-1)) }^F$
is
the $F$-distribution with $(a-1,ab(n-1)) $ degrees of freedom.
Thus, as seen in the formula
(\ref{eq46}),
it suffices to calculate
the $\alpha$-point
$F_{ab(n-1), \alpha}^{a-1}$
Thus,
we see
\begin{align}
(\eta^\alpha_{\omega})^2 =
F_{ab(n-1), \alpha}^{a-1}
\cdot
n(a-1)/(n-1)
\label{eq67}
\end{align}
\par
\noindent
Therefore,
we get ${\widehat R}_{\widehat{x}}^{\alpha; \Theta}$
(or,
${\widehat R}_{\widehat{x}}^{\alpha; X}$;
the $({}\alpha{})$-rejection region
of
$H_N =\{(0.0. \ldots, 0)\}( \subseteq \Theta= {\mathbb R}^a)$
)
as follows:
\begin{align}
{\widehat R}_{{H_N}}^{\alpha; \Theta}
&
=
\bigcap_{\omega =((\mu_i)_{i=1}^a, \sigma ) \in  \Omega (={\mathbb R}^a \times {\mathbb R}_+ ) \mbox{ \footnotesize such that }
\pi(\omega)= (\alpha_i)_{i=1}^a \in {H_N}=\{(0,0,\ldots,0)\}}
\{
E({x})
(\in
\Theta)
:
d_\Theta^{x}  ({}E({x}),
\pi(\omega))
\ge
\eta^\alpha_{\omega }
\}
\nonumber
\\
&
=
\{
E({x})
(\in
\Theta)
:
\frac{
(\sum_{i=1}^a \sum_{j=1}^b(
x_{ij \cdot}
- x_{\cdot \cdot \cdot}
)^2)/(a-1)}{
(\sum_{i=1}^a \sum_{j=1}^b\sum_{k=1}^n (x_{ijk} - x_{ij \cdot})^2)
/(ab(n-1))
}
\ge
F_{ab(n-1), \alpha}^{a-1}
\}
\label{eq68}
\end{align}
Thus,
\begin{align}
{\widehat R}_{{H_N}}^{\alpha; X}=
E^{-1}({\widehat R}_{{H_N}}^{\alpha; \Theta})
=
\{
x
(\in
X)
:
\frac{
(\sum_{i=1}^a \sum_{j=1}^b(
x_{ij \cdot}
- x_{\cdot \cdot \cdot}
)^2)/(a-1)}{
(\sum_{i=1}^a \sum_{j=1}^b\sum_{k=1}^n (x_{ijk} - x_{ij \cdot})^2)
/(ab(n-1))
}
\ge
F_{ab(n-1), \alpha}^{a-1}
\}
\label{eq69}
\end{align}
Therefore, the statistical hypothesis test (D$_2$) in Theorem 1 is applicable.

\rm
\vskip0.5cm
\par

\par
\noindent
\bf
Remark 3
\rm
If we assume the null hypothesis such that $\beta_1=\beta_2=\cdots=\beta_b=0$,
we can give the similar answer such as Example 4.
\rm
\par
\noindent
\vskip0.5cm
\par
\noindent
\bf
Example 5
\rm
[The null hypothesis such that $(\alpha \beta)_{ij}=0$
($\forall i=1,2, \ldots, a,\;j=1,2, \ldots, b$)].
\rm
Put
\begin{align}
\Theta = {\mathbb R}^{ab}
\label{eq70}
\end{align}
and define the $\pi : \Omega \to \Theta $
such that
\begin{align}
\Omega = {\mathbb R}^a \times {\mathbb R}_+
\ni
\omega =(\mu_1, \mu_2, \ldots, \mu_a, \sigma)
\mapsto
\pi(\omega)
=
((\alpha \beta)_{ij})_{i=1,2, \ldots, a,\;\;j=1,2, \ldots, b }
\in
\Theta = {\mathbb R}^{ab}
\label{eq71}
\end{align}
where, as defined in (\ref{eq50}),
\begin{align}
(\alpha \beta)_{ij}
=\mu_{ij}-\mu_{i \cdot}
-\mu_{\cdot j }+\mu_{\cdot \cdot }
\label{eq72}
\end{align}
Define the null hypothesis
$H_N ( \subseteq \Theta = {\mathbb R}^{ab})$
such that
\begin{align}
H_N
&
=
\{
((\alpha \beta)_{ij})_{i=1,2, \ldots, a,\;\;j=1,2, \ldots, b }
\in
\Theta = {\mathbb R}^{ab}
\;:\;
(\alpha \beta)_{ij}=0,
(\forall {i=1,2, \ldots, a,\;\;j=1,2, \ldots, b }
)
\}
\label{eq73}
\end{align}
Put
\begin{align}
&
\|
\theta^{(1)}- \theta^{(2)}
\|_\Theta
=
\sqrt{
\sum_{i=1}^a \sum_{j=1}^b \Big(\theta_{ij}^{(\ell)} - \theta_{ij}^{(\ell)} 
\Big)^2
}
\label{eq74}
\\
&
\qquad
(\forall \theta^{(\ell)} =( \theta_{ij}^{(\ell)})_{i=1,2, \ldots, a,\;\;j=1,2, \ldots, b }
\in {\mathbb R}^{ab},
\;
\ell=1,2
)
\nonumber
\end{align}

\par
\noindent
Define 
$
{\overline{SS}}(x)
$
by the formula (\ref{eq58}),
and
define
the semi-distance
$d_\Theta^x$ in $\Theta$ such that
\begin{align}
&
d_\Theta^x (\theta^{(1)}, \theta^{(2)})
=
\frac{\|\theta^{(1)}- \theta^{(2)} \|_\Theta}{
\sqrt{{\overline{SS}}(x)}
}
\qquad
(\forall \theta^{(1)}, \theta^{(2)}
\in \Theta,
 \forall x \in X )
 \label{eq75}
\end{align}

\par
\noindent
Define and calculate the estimator
$E: X(={\mathbb R}^{abn}) \to \Theta(={\mathbb R}^{ab} )$ such that
\begin{align}
&
E(
(x_{ijk})_{i=1,...,a, \; j=1,2,...b, \; k=1,2,...,n })
\nonumber
\\
=
&
\Big(\frac{
\sum_{k=1}^n x_{ijk}}{n}
-
\frac{
\sum_{j=1}^b \sum_{k=1}^n x_{ijk}}{bn}
-
\frac{
\sum_{j=1}^b\sum_{k=1}^n x_{ijk}}{an}
+
\frac{\sum_{i=1}^a \sum_{j=1}^b \sum_{k=1}^n x_{ijk}}{abn}
\Big)_{i=1,2,\ldots, a \; j=1,2,...b, }
\nonumber
\\
=
&
\Big( x_{i j \cdot }
-
x_{i \cdot \cdot }
-
x_{\cdot j \cdot }
+
x_{\cdot \cdot \cdot }
\Big)_{i=1,2,\ldots, a \; j=1,2,...b,}
\label{eq76}
\intertext{and thus,}
&
E(
(x_{ijk}-\mu_{ij})_{i=1,...,a, \; j=1,2,...b, \; k=1,2,...,n })
\nonumber
\\
=
&
\Big(\frac{
\sum_{k=1}^n (x_{ijk}-\mu_{ij})}{n}
-
\frac{
\sum_{j=1}^b \sum_{k=1}^n (x_{ijk}-\mu_{ij})}{bn}
\nonumber
\\
&
\qquad
-
\frac{
\sum_{j=1}^b\sum_{k=1}^n (x_{ijk}-\mu_{ij})}{an}
+
\frac{\sum_{i=1}^a \sum_{j=1}^b \sum_{k=1}^n (x_{ijk}-\mu_{ij})}{abn}
\Big)_{i=1,2,\ldots, a \; j=1,2,...b, }
\nonumber
\\
=
&
\Big( (x_{i j} \cdot - \mu_{ij})
-
(x_{i \cdot \cdot }-\mu_{i \cdot})
-
(x_{\cdot j \cdot }- \mu_{\cdot j})
+
(x_{\cdot \cdot \cdot }- \mu_{\cdot \cdot })
\Big)_{i=1,2,\ldots, a \; j=1,2,...b,}
\nonumber
\\
=
&
\Big( x_{i j \cdot }
-
x_{i \cdot \cdot }
-
x_{\cdot j \cdot }
+
x_{\cdot \cdot \cdot }
\Big)_{i=1,2,\ldots, a \; j=1,2,...b}
\qquad
(\mbox{under the null hypothesis }(\alpha \beta)_{ij}=0 )
\label{eq77}
\end{align}

Therefore,
\begin{align}
E(
(x_{ijk})_{i=1,...,a, \; j=1,2,...b, \; k=1,2,...,n })
=
E(
(x_{ijk}-\mu_{ij})_{i=1,...,a, \; j=1,2,...b, \; k=1,2,...,n })
\label{eq78}
\end{align}

\par

\par
\noindent
Hence, for each
$i=1,...,a, \; j=1,2,...b,$
\begin{align}
&
E_{ij}
(x_{ijk} -\mu_{ij})
\nonumber
\\
=
&
\frac{
\sum_{k=1}^n (x_{ijk}-\mu_{ij})}{n}
-
\frac{
\sum_{j=1}^b \sum_{k=1}^n (x_{ijk}-\mu_{ij})}{bn}
-
\frac{
\sum_{j=1}^b\sum_{k=1}^n (x_{ijk}-\mu_{ij})}{an}
\nonumber
\\
&
\qquad \qquad
+
\frac{\sum_{i=1}^a \sum_{j=1}^b \sum_{k=1}^n(x_{ijk}-\mu_{ij})}{abn}
\nonumber
\\
=
&
E_{ij}(x)- (\alpha\beta)_{ij}
\nonumber
\\
=
&
x_{i j \cdot }
-
x_{i \cdot \cdot }
-
x_{\cdot j \cdot }
+
x_{\cdot \cdot \cdot }
-
(\alpha\beta)_{ij}
\label{eq79}
\end{align}

\par
\noindent
Thus, we see that
\begin{align}
&
\| E(x) - \pi (\omega )\|^2_\Theta
\nonumber
\\
=
&
||
\Big(
E_{ij}(x) - (\alpha \beta)_{ij}
\Big)_{i=1,2,\ldots, a \; j=1,2,...b }
||_\Theta^2
\label{eq80}
\intertext{and thus, if the null hypothesis $H_N$ is assumed
(i.e.,
$(\alpha \beta)_{ij}=0$
$(\forall i=1,2,\ldots, a, \; j=1,2, \ldots, b )$
)
}
=
&
\sum_{i=1}^a\sum_{j=1}^b(
x_{i j \cdot }
-
x_{i \cdot \cdot }
-
x_{\cdot j \cdot }
+
x_{\cdot \cdot \cdot }
)^2
\label{eq81}
\end{align}

\par
\noindent
Thus,
for any
$ \omega=(\mu_1, \mu_2 )  ({}\in \Omega=
{\mathbb  R} \times {\mathbb R} )$,
define
the positive number
$\eta^\alpha_{\omega}$
$({}> 0)$
such that:
\begin{align}
\eta^\alpha_{\omega}
=
\inf
\{
\eta > 0:
[G({}E^{-1} ({}
{{\rm Ball}^C_{d_\Theta^{x}}}(\pi(\omega) ; \eta{}))](\omega )
\ge \alpha
\}
\label{eq82}
\end{align}

\par
\noindent
Assume the null hypothesis $H_N$
(i.e.,
$(\alpha \beta)_{ij}=0$
$(\forall i=1,2,\ldots, a, \; j=1,2, \ldots, b )$
).
Now let us calculate the $\eta^\alpha_{\omega}$ as follows:
\begin{align}
&
E^{-1}({{\rm Ball}^C_{d_\Theta^{x} }}(\pi(\omega) ; \eta ))
=\{ x \in X = {\mathbb R}^{abn}
\;:\;
d_\Theta^x (E(x), \pi(\omega ))
>
\eta
\}
\nonumber
\\
=
&
\{ x \in X = {\mathbb R}^{abn}
\;:\;
\frac{
abn \sum_{i=1}^a\sum_{j=1}^b(
x_{i j \cdot }
-
x_{i \cdot \cdot }
-
x_{\cdot j \cdot }
+
x_{\cdot \cdot \cdot }
)^2
}{
\sum_{i=1}^a \sum_{j=1}^b\sum_{k=1}^n (x_{ijk} - x_{ij \cdot})^2
}
>
\eta^2
\}
\label{eq83}
\end{align}

\par
\noindent
That is,
for any $\omega
=((\mu_{ij})_{i=1,2,\ldots,a, \;j=1,2,\ldots,b},\;, \sigma) \in \Omega
={\mathbb R}^{ab +1}$
such tht
$\pi( \omega )\in H_N
(\subseteq
{\mathbb R}^{ab}
)
$
(i.e.,
$(\alpha \beta)_{ij}=0$
$(\forall i=1,2,\ldots, a, \; j=1,2, \ldots, b )$
)
\par
\noindent
\begin{align}
&
[{{{G}}}^{abn}
(
E^{-1}({{\rm Ball}^C_{d_\Theta^{x} }}(\pi(\omega) ; \eta ))
)
({}\omega{})
\nonumber
\\
=
&
\frac{1}
{({
{\sqrt{2 \pi }
\sigma}
})^{abn}}
\underset{
E^{-1}({{\rm Ball}^C_{d_\Theta^{x} }}(\pi(\omega) ; \eta ))
}
{\int \cdots \int}
\exp[- 
\frac{
\sum_{i=1}^a \sum_{j=1}^b \sum_{k=1}^n  (x_{ijk} - 
\mu_{ij}
)^2 
}{2 \sigma^2}
] 
\bigtimes_{k=1}^n
\bigtimes_{j=1}^b
\bigtimes_{i=1}^a
d{x_{ijk} }
\nonumber
\\
=
&
\frac{1}
{({
{\sqrt{2 \pi }
\sigma}
})^{abn}}
\underset{
\{x \in X \;:\; d_\Theta^x ( E(x), \pi(\omega ) \ge \eta \}
}
{\int \cdots \int}
\exp[- 
\frac{
\sum_{i=1}^a \sum_{j=1}^b \sum_{k=1}^n  (x_{ijk} - 
\mu_{ij}
)^2 
}{2 \sigma^2}
] 
\bigtimes_{k=1}^n
\bigtimes_{j=1}^b
\bigtimes_{i=1}^a
d{x_{ijk} }
\nonumber
\\
=
&
\frac{1}
{({
{\sqrt{2 \pi }}
})^{abn}}
\underset{
\frac{
\sum_{i=1}^a\sum_{j=1}^b(
x_{i j \cdot }
-
x_{i \cdot \cdot }
-
x_{\cdot j \cdot }
+
x_{\cdot \cdot \cdot }
)^2
}{
\sum_{i=1}^a \sum_{j=1}^b\sum_{k=1}^n (x_{ijk} - x_{ij \cdot})^2
}
>
\frac{\eta^2}{abn}
}
{\int \cdots \int}
\exp[- 
\frac{
\sum_{i=1}^a \sum_{j=1}^b \sum_{k=1}^n  (x_{ijk} 
)^2 
}{2}
] 
\bigtimes_{k=1}^n
\bigtimes_{j=1}^b
\bigtimes_{i=1}^a
d{x_{ijk} }
\nonumber
\\
=
&
\frac{1}
{({
{\sqrt{2 \pi }}
})^{abn}}
\underset{
\frac{
\frac{
\sum_{i=1}^a\sum_{j=1}^b(
x_{i j \cdot }
-
x_{i \cdot \cdot }
-
x_{\cdot j \cdot }
+
x_{\cdot \cdot \cdot }
)^2}{(a-1)(b-1)}
}{
\frac{\sum_{i=1}^a \sum_{j=1}^b\sum_{k=1}^n (x_{ijk} - x_{ij \cdot})^2}{
ab(n-1)
}
}
>
\frac{\eta^2(ab(n-1))}{abn(a-1)(b-1)}
}
{\int \cdots \int}
\exp[- 
\frac{
\sum_{i=1}^a \sum_{j=1}^b \sum_{k=1}^n  (x_{ijk} 
)^2 
}{2}
] 
\bigtimes_{k=1}^n
\bigtimes_{j=1}^b
\bigtimes_{i=1}^a
d{x_{ijk} }
\label{eq84}
\intertext{
\begin{itemize}
\item[(E$_4$)]
using the formula of Gauss integrals derived in Kolmogorov's probability theory
(also, recall Remark 2), we finally get as follows.
\end{itemize}} 
=
&
\int^{\infty}_{\frac{\eta^2 (n-1)}{n (a-1)(b-1)}} 
p_{((a-1)(b-1),ab(n-1)) }^F
(t) dt =\alpha
(\mbox{  e.g., } \alpha=0.05)
\label{eq85}
\end{align}
where
$p_{((a-1)(b-1),ab(n-1)) }^F$
is
the $F$-distribution with $((a-1)(b-1),ab(n-1)) $ degrees of freedom.
Thus, as seen in the formula
(\ref{eq67}),
Thus, it suffices to put
\begin{align}
{\frac{\eta^2 (n-1)}{n (a-1)(b-1)}}
={F_{ab(n-1), \alpha}^{(a-1)(b-1)} } 
(=\mbox{"$\alpha$-point"})
\label{eq86}
\end{align}
And thus we see,
\begin{align}
(\eta^\alpha_{\omega})^2 =
{F_{ab(n-1), \alpha}^{(a-1)(b-1)} } 
n(a-1)(b-1)/(n-1)
\label{eq87}
\end{align}
\par
\noindent
%
%
%
\par
\noindent
Therefore,
we get ${\widehat R}_{\widehat{x}}^{\alpha; \Theta}$
(or,
${\widehat R}_{\widehat{x}}^{\alpha; X}$;
the $({}\alpha{})$-rejection region
of
$H_N =\{((\alpha \beta)_{ij})_{i=1,2, \cdots, a,
j=1,2, \cdots, b}
\; :\;
(\alpha \beta)_{ij}=0
\;
(i=1,2, \cdots, a,
j=1,2, \cdots, b)\}( \subseteq \Theta= {\mathbb R}^{ab})$
)
as follows:
\begin{align}
{\widehat R}_{{H_N}}^{\alpha; \Theta}
&
=
\bigcap_{\omega =((\mu_{ij})_{i=1}^a{}_{j=1}^b, \sigma ) \in  \Omega (={\mathbb R}^a \times {\mathbb R}_+ ) \mbox{ \footnotesize such that }
\pi(\omega)= (\alpha \beta)_{ij} \in {H_N}}
\{
E({x})
(\in
\Theta)
:
d_\Theta^{x}  ({}E({x}),
\pi(\omega))
\ge
\eta^\alpha_{\omega }
\}
\nonumber
\\
&
=
\{
E({x})
(\in
\Theta)
:
\frac{
(\sum_{i=1}^a \sum_{j=1}^b(
x_{ij \cdot}
- x_{\cdot \cdot \cdot}
)^2)/((a-1)(b-1))}{
(\sum_{i=1}^a \sum_{j=1}^b\sum_{k=1}^n (x_{ijk} - x_{ij \cdot})^2)
/(ab(n-1))
}
\ge
{F_{ab(n-1), \alpha}^{(a-1)(b-1)} } 
\}
\label{eq88}
\end{align}
Thus,
\begin{align}
{\widehat R}_{{H_N}}^{\alpha; X}=
E^{-1}({\widehat R}_{{H_N}}^{\alpha; \Theta})
=
\{
x
(\in
X)
:
\frac{
(\sum_{i=1}^a \sum_{j=1}^b(
x_{ij \cdot}
- x_{\cdot \cdot \cdot}
)^2)/((a-1)(b-1))}{
(\sum_{i=1}^a \sum_{j=1}^b\sum_{k=1}^n (x_{ijk} - x_{ij \cdot})^2)
/(ab(n-1))
}
\ge
{F_{ab(n-1), \alpha}^{(a-1)(b-1)} } 
\}
\label{eq89}
\end{align}
Therefore, the statistical hypothesis test (D$_2$) in Theorem 1 is applicable.

\rm
\vskip0.5cm
\par

\section{Conclusions
}
\par
\noindent
\par

\par
\noindent
\par



We believe that quantum language has a great powet of
description,
and therefore, even statistics can be described by quantum language.
%
%
%
Since quantum language is suited for theoretical arguments, we believe, from the theoretical point of view, that our results
(i.e.,
ANOVA in Section 2) are visible and simple.
Therefore,
we can easily answer the following question:
\begin{itemize}
\item[(F$_1$)]
Where is Kolmogorov's probability theory used in ANOVA?
\end{itemize}
As the conclusion, we can answer as follows:
\begin{itemize}
\item[(F$_2$)]
Kolmogorov's probability theory
is merely used in order to calculate
multi-dimensional Gauss integrals
throughout this paper
({\it cf}.
the items
(E$_1$)-(E$_4$)
in Examples 2-5).
\end{itemize}
It is reasonable, since Kolmogorov's probability theory is mathematics.
Although we may calculate the multi-dimensional Gauss integrals
without Kolmogorov's probability theory ({\it cf.} Remark 2),
it is sure that
the conventional calculation
(due to Kolmogorov's probability theory)
is elegant and powerful.
In this sense,
we believe that
mathematical theories
(particularly,
Kolmogorov's probability theory
and the theory of operator algebra
({\it cf}. \cite{Saka})
)
are
indispensable for quantum language.
\par
We hope that our assertions will be examined from various points of view.

%
%



\rm
\par
\renewcommand{\refname}{
\large 
References}
{
\small

\normalsize
}

\end{document}